\documentclass[letterpaper, 10 pt, conference]{ieeeconf}  

\IEEEoverridecommandlockouts                              
\usepackage{graphics} 
\usepackage{amsmath} 
\usepackage{amssymb}  

\usepackage{cite}
\usepackage{graphicx}
\graphicspath{{graphics/}}
\usepackage{subfigure}
\usepackage{enumerate}
\usepackage{ifpdf}
\usepackage{tikz}
\usepackage{tikz-qtree}
\usetikzlibrary{arrows,petri,topaths,automata}
\usepackage{tikz,fullpage}
\usepackage{pgf}
\usepackage{color}
\usepackage{algorithmic} 
\usepackage{algorithm} 
\usepackage[textsize=tiny,textwidth=2cm]{todonotes}	

\newcommand{\be}{\begin{enumerate}}
\newcommand{\ee}{\end{enumerate}}

\newcommand{\Tr}{\text{Tr}}
\newcommand{\TF}[4]{\left [ \begin{array}{c|c}{#1} & {#2} \\ \hline {#3} & {#4} \end{array} \right ]}

\newcommand{\D}{\mathsf D}
\newcommand{\R}{\mathsf R}
\newcommand{\Adj}{\text{Adj}}

\newcommand{\Prob}{\mathbb P}
\newcommand{\Exp}{\mathbb E}

\newtheorem{thm}{Theorem}

\newtheorem{cor}{Corollary}
\newtheorem{definition}{Definition}

\newtheorem{exmp}{Example}[section]
\newtheorem{rem}{Remark}
\newtheorem{remark}{Remark}

\title{\LARGE \bf
Differentially Private Kalman Filtering
}

\author{Jerome Le Ny and George J. Pappas
\thanks{J. Le Ny is with the department of Electrical Engineering, Ecole Polytechnique de Montreal,
QC H3C 3A7, Canada. G.~Pappas is with the Department of Electrical and Systems Engineering, 
University of Pennsylvania, Philadelphia, PA 19104, USA. {\tt\small jerome.le-ny@polymtl.ca, 
pappasg@seas.upenn.edu.}}%
}

\begin{document}

\maketitle
\thispagestyle{empty}
\pagestyle{empty}

\begin{abstract}
This paper studies the $\mathcal H_2$ (Kalman) filtering problem in the situation
where a signal estimate must be constructed based on inputs from individual participants, 
whose data must remain private. This problem arises in emerging applications such as smart grids 
or intelligent transportation systems, where users continuously send data to third-party 
aggregators performing global monitoring or control tasks, and require guarantees 
that this data cannot be used to infer additional personal information. 
To provide strong formal privacy guarantees against adversaries with arbitrary side information, 
we rely on the notion of \emph{differential privacy} introduced relatively recently in the database 
literature.
This notion is extended to dynamic systems with many participants contributing independent 
input signals, and mechanisms are then proposed to solve the $\mathcal H_2$ filtering problem 
with a differential privacy constraint.
A method for mitigating the impact of the privacy-inducing mechanism on the estimation performance 
is described, which relies on controlling the $\mathcal H_\infty$ norm of the filter. 
Finally, we discuss an application to a privacy-preserving traffic monitoring system.
\end{abstract}





\section{Introduction}

In many applications, such as smart grids, population health monitoring, or traffic monitoring, 
the efficiency of the system relies on the participation of the users to provide reliable data 
in real-time, e.g., power consumption, sickness symptoms, or GPS coordinates.
However, for privacy or security reasons, the participants benefiting from these services 
generally do not want to release more information than strictly necessary. 
Unfortunately, examples of unintended loss of privacy already abound. Indeed, it is possible 
to infer from the trace of a smart meter the type of appliances present in a house as well 
as the occupants' daily activities \cite{Hart92_loadMonitoring}, to re-identify an anonymous 
GPS trace by correlating it with publicly available information such as location 
of work \cite{Hoh11_VTL_trafficMonitoring}, or to infer individual transactions on commercial 
websites from temporal changes in public recommendation systems \cite{Calandrino11_privacyAttackCollabFilt}.
Providing rigorous guarantees to the users about the privacy risks incurred is thus crucial 
to encourage participation and ultimately realize the benefits promised by these systems.

In our recent work \cite{LeNy_CDC12_DP}, we introduced privacy concerns in the the context 
of systems theory, by relying on the notion of \emph{differential privacy} \cite{Dwork06_DPcalibration}, 
a particularly successful definition of privacy used in the database literature.
This notion is motivated by the fact that any useful information provided by a dataset about a group 
of people can compromise the privacy of specific individuals due to the existence of side information. 
Differentially private mechanisms randomize their responses to dataset analysis requests and 
guarantee that whether or not an individual chooses to contribute her data only marginally 
changes the distribution over the published outputs. As a result, even an adversary cross-correlating
these outputs with other sources of information cannot infer much more about specific individuals after 
publication than before \cite{Kasiviswanathan08_DP_sideInfo}.

Most work related to privacy is concerned with the analysis of static databases, whereas cyber-physical 
systems clearly emphasize the need for mechanisms working with dynamic, time-varying data streams. 
Recently, information-theoretic approaches have been proposed to guarantee some level of privacy
when releasing time series \cite{Varodayan11_batteryPrivacy, Sankar11_privacyInfoTheoretic}. 
However, the resulting privacy guarantees only hold if the statistics of the participants' data streams 
obey the assumptions made (typically stationarity, dependence and distributional assumptions), 
and require the explicit statistical modeling of all available side information. This task is impossible 
in general as new, as-yet-unknown side information can become available after releasing the results. 
In contrast, differential privacy is a worst-case notion that holds independently of any probabilistic 
assumption on the dataset, and controls the information leakage against adversaries with arbitrary 
side information \cite{Kasiviswanathan08_DP_sideInfo}. Once such a privacy guarantee is enforced, 
one can still leverage potential additional statistical information about the dataset to improve the quality 
of the outputs. 

In this paper, we pursue our work on differential privacy for dynamical systems \cite{LeNy_CDC12_DP}, 
by considering the $\mathcal H_2$ filtering problem (or steady-state Kalman filtering) with a differential privacy constraint. 
In this problem, the goal is to minimize an estimation error variance for a desired linear combination of the participants' state 
trajectories, based on their contributed measurements, while guaranteeing the privacy of the individual signals.
In contrast to the generic filtering mechanisms presented in \cite{LeNy_CDC12_DP}, we emphasize here
how a model of the participants' dynamics can be leveraged to publish more accurate results, 
without compromising the differential privacy guarantee if this model is not accurate.
Section \ref{section: differential privacy background} provides some technical background on differential privacy
and Section \ref{section: DP linear systems} describes a basic mechanism enforcing privacy for dynamical systems
by injecting additional white noise.
As shown in \cite{LeNy_CDC12_DP}, accurate private results can be published for filters with small incremental gains 
with respect to the individual input channels. This leads us in Section \ref{section: private KF} to present 
a modification of the standard Kalman filter, essentially controlling its $\mathcal H_\infty$ norm simultaneously 
with the steady-state estimation error, in order to minimize the impact of the privacy-inducing mechanism.
Finally, Section \ref{section: traffic monitoring example} describes an application to a simplified 
traffic monitoring system relying on location traces from the participants to provide an average velocity
estimate on a road segment. Most proofs are omitted from this extended abstract and will appear in the
full version of the paper.


\section{Differential Privacy}	\label{section: differential privacy background}

In this section we review the notion of differential privacy \cite{Dwork06_DPcalibration} 
as well as a basic mechanism that can be used to achieve it when the released data 
belongs to a finite-dimensional vector space. We refer the reader to the surveys by Dwork, 
e.g., \cite{Dwork_ICAL06_DP}, for additional background on differential privacy.


\subsection{Definition}

Let us fix some probability space $(\Omega, \mathcal F, \Prob)$. 
Let $\D$ be a space of datasets of interest 
(e.g., a space of data tables, or a signal space).
A \emph{mechanism} is just a map $M: \D \times \Omega \to \R$, for some
measurable output space $\R$, such that for any element $d \in \D$, 
$M(d,\cdot)$ is a random variable, typically writen simply $M(d)$. 
A mechanism can be viewed as a probabilistic algorithm to 
answer a query $q$, which is a map $q: \D \to \R$. 
In some cases, we index the mechanism by the query $q$ of interest, writing $M_q$.

\begin{exmp}	\label{ex: space of databases}
Let $\D = \mathbb R^{n}$, with each real-valued entry of $d \in \D$ corresponding to some sensitive
information for an individual contributing her data.
A data analyst would like to know the average of the entries of $d$, i.e., her query is
\begin{align}
&q: \D \to \mathbb R \nonumber, \;\; q(d) = \frac{1}{n} \sum_{i=1}^n d_{i}. \nonumber
\end{align}
As detailed in Section \ref{section: basic mech}, a typical 
mechanism $M_q$ to answer this query in a 
differentially private way computes $q(d)$ and blurs 
the result by adding a random variable $Y: \Omega \to \mathbb R$
\begin{align}
&M_q: \D \times \Omega \to \mathbb R, 
\;\; M_{q}(d) = \frac{1}{n} \sum_{i=1}^n d_{i} + Y.	\nonumber 
\end{align}
Note that in the absence of perturbation $Y$, an adversary who knows $n$ and $d_j, j \geq 2,$
can recover the remaining entry $d_1$ exactly if he learns $q(d)$. This can deter people
from contributing their data, even though broader participation improves the accuracy 
of the analysis and thus can be beneficial to the population as a whole.
\end{exmp}

Next, we introduce the definition of differential privacy.  We call a measure $\mu$
on $\R$ $\delta$-bounded if it is a finite positive measure with $\mu(\R) \leq \delta$.
Intuitively in the following definition, $\D$ is a space of datasets of interest,
and we have a binary relation $\Adj$ on $\D$, called adjacency, such that $\Adj(d,d')$ if and only if 
$d$ and $d'$ differ by the data of a single participant. 

\begin{definition}	\label{def: differential privacy original}
Let $\D$ be a space equipped with a binary relation denoted $\Adj$,
and let $(\R, \mathcal M)$ be a measurable space.
Let $\epsilon, \delta \geq 0$. 
A mechanism $M: \D \times \Omega \to \R$ is 
$(\epsilon, \delta)$-differentially private if 
there exists a $\delta$-bounded measure $\mu$ on $(\R,\mathcal M)$ such that
for all $d,d' \in \D$ such that $\Adj(d,d')$ and for all $S \in \mathcal M$, we have
\begin{align}	\label{eq: standard def approximate DP original}
\Prob(M(d) \in S) \leq e^{\epsilon} \Prob(M(d') \in S) + \mu(S). 
\end{align}
If $\delta=0$, the mechanism is said to be $\epsilon$-differentially private. 
\end{definition}

This definition is essentially the same as the one introduced in \cite{Dwork06_DPcalibration} 
and subsequent work, except for the fact that $\mu(S)$ in (\ref{eq: standard def approximate DP original}) 
is usually replaced by the constant $\delta$. The definition says that for two adjacent datasets, 
the distributions over the outputs of the mechanism should be close.  The choice of the parameters $\epsilon, \delta$ is set by the privacy policy.
Typically $\epsilon$ is taken to be a small constant, e.g., $\epsilon \approx 0.1$ or perhaps even $\ln 2$ or $\ln 3$.
The parameter $\delta$ should be kept small as it controls the probability of
certain significant losses of privacy, e.g., when a zero probability event for input $d'$ becomes an event
with positive probability for input $d$ in (\ref{eq: standard def approximate DP original}).
%

\begin{remark}
The definition of differential privacy depends on the choice of $\sigma$-algebra $\mathcal M$ 
in Definition \ref{def: differential privacy original}. When we need to state this $\sigma$-algebra
explicitly, we write $M: \D \times \Omega \to (\R,\mathcal M)$.
In particular, this $\sigma$-algebra should be
sufficiently ``rich'', since (\ref{eq: standard def approximate DP original}) is trivially satisfied
by any mechanism if $\mathcal M = \{\emptyset, \R\}$. 
\end{remark}


A useful property of the notion of differential privacy is that 
no additional privacy loss can occur by simply manipulating 
an output that is differentially private. This result is similar 
in spirit to the data processing inequality from information theory \cite{Cover91_infoTheory}.

\begin{thm}[resilience to post-processing]	\label{thm: resilience to post-processing}
Let $M_1: \D \times \Omega \to (\R_1,\mathcal M_1)$ be an $(\epsilon,\delta)$-differentially private mechanism.
Let $M_2: \D \times \Omega \to (\R_2,\mathcal M_2)$ be another mechanism 
such that for all $S \in \mathcal M_2$, there exists a nonnegative 
measurable function $f_S$ such that for all $d \in \D$, we have
\begin{align}	\label{eq: post-processing definition}
\Prob(M_2(d) \in S | M_1(d)) = f_S(M_1(d)), \forall d \in \D. 
\end{align}
Then $M_2$ is $(\epsilon,\delta)$-differentially private.
\end{thm}

\vspace{0.2cm}

\begin{rem}
Suppose that $M_1$ takes its values in a discrete set. Then the condition (\ref{eq: post-processing definition})
says that the conditional distribution $\Prob(M_2(d) \in S | M_1(d) = m_1)$ for a given element $m_1$ 
does not further depend of $d$. In other words, a mechanism $M_2$ accessing the dataset 
only indirectly via the output of $M_1$ cannot weaken the privacy guarantee. 
Hence post-processing can be used to improve the \emph{accuracy} of an output, 
without weakening the privacy guarantee.
\end{rem}



\subsection{A Basic Differentially Private Mechanism}	\label{section: basic mech}

A mechanism that throws away all the information in a dataset is obviously private, 
but not useful, and in general one has to trade off privacy for utility when 
answering specific queries. 
We recall below a basic mechanism that can be used to
answer queries in a differentially private way.
%
%
We are only concerned in this section with queries that return numerical answers, 
i.e., here a query is a map $q: \D \to \R$, where the output space $\R$ equals $\mathbb R^k$ for some $k > 0$,
is equipped with a norm  denoted $\| \cdot \|_\R$,
and the $\sigma$-algebra $\mathcal M$ on $\R$ is taken to be the standard 
Borel $\sigma$-algebra, denoted $\mathcal R^k$. 
The following quantity plays an important
role in the design of differentially private mechanisms \cite{Dwork06_DPcalibration}.

\begin{definition}	\label{defn: sensitivity}
Let $\D$ be a space equipped with an adjacency relation $\Adj$.
The sensitivity of a query $q: \D \to \R$ is defined as
\[
\Delta_\R q := \max_{d,d':\Adj(d,d')} \|q(d) - q(d')\|_\R.
\]
In particular, for $\R = \mathbb R^k$ equipped with the $p$-norm 
$\| x \|_p = \left(\sum_{i=1}^k |x_i|^p \right)^{1/p}$, 
for $p \in [1,\infty]$,
we denote the $\ell_p$ sensitivity by $\Delta_p q$.
\end{definition}
A differentially private mechanism proposed in  \cite{Dwork06_DPgaussian}, modifies an answer to a 
numerical query by adding iid zero-mean noise distributed according to a Gaussian distribution.
Recall the definition of the $\mathcal Q$-function 
\[
\mathcal Q(x) := \frac{1}{\sqrt{2 \pi}} \int_x^{\infty} e^{-\frac{u^2}{2}} du.
\]
The following theorem tightens the analysis from \cite{Dwork06_DPgaussian}.

\begin{thm}	\label{thm: Gaussian mech}
Let $q: \D \to \mathbb R^k$ be a query.
Then the Gaussian mechanism $M_Q: \D \times \Omega \to \mathbb R^k$ 
defined by $M_q(d) = q(d) + w$, with $w \sim \mathcal N\left(0,\sigma^2 I_k \right)$, 
where $\sigma \geq \frac{\Delta_2 q}{2 \epsilon}(K + \sqrt{K^2+2\epsilon})$ and $K = \mathcal Q^{-1}(\delta)$,
is $(\epsilon,\delta)$-differentially private.
\end{thm}

For the rest of the paper, we define 
\[
\kappa(\delta,\epsilon) = \frac{1}{2 \epsilon} (K+\sqrt{K^2+2\epsilon}),
\]
so that the standard deviation $\sigma$ in Theorem \ref{thm: Gaussian mech}
can be written $\sigma(\delta,\epsilon) = \kappa(\epsilon,\delta) \Delta_2 q$.
It can be shown that $\kappa(\delta,\epsilon)$ behaves roughly as $O(\ln(1/\delta))^{1/2}/\epsilon$.
For example, to guarantee $(\epsilon,\delta)$-differential privacy with $\epsilon = \ln(2)$ and $\delta = 0.05$, 
we obtain that the standard deviation of the Gaussian noise introduced should be about $2.65$ times 
the $\ell_2$-sensitivity of $q$.


\section{Differentially Private Dynamic Systems}		\label{section: DP linear systems}

In this section we review the notion of differential privacy for dynamic systems, 
following \cite{LeNy_CDC12_DP}.
We start with some notations and technical prerequisites.
All signals are discrete-time signals and all systems are assumed to be causal.
For each time $T$, let $P_T$ be the truncation operator, so that for any signal $x$ we have
\[
(P_T x)_t = \begin{cases}
x_t, & t \leq T \\
0, & t > T.
\end{cases}
\]
Hence a deterministic system $\mathcal G$ is causal if and only if $P_T \mathcal G = P_T \mathcal G P_T$.
We denote by $\ell_{p,e}^m$ the space of sequences with values in
$\mathbb R^m$ and such that $x \in \ell_{p,e}^m$ if and only if $P_T x$
has finite $p$-norm for all integers $T$. The $\mathcal H_2$ norm and $\mathcal H_\infty$
norm of a stable transfer function $\mathcal G$ are defined respectively as
\begin{align*}
\|\mathcal G\|_2 &= \left(\frac{1}{2 \pi} \int_{-\pi}^\pi \Tr (\mathcal G^*(e^{i \omega}) \mathcal G(e^{i \omega})) d \omega\right)^{1/2}, \\
\|\mathcal G\|_\infty &= \text{ess} \hspace{-0.3cm} \sup_{\omega \in [-\pi,\pi)} \sigma_{\max} (\mathcal G(e^{i \omega})),
\end{align*}
where $\sigma_{\max}(A)$ denotes the maximum singular value of a matrix $A$.

We consider situations in which private participants contribute input signals 
driving a dynamic system and the queries consist of output signals 
of this system. 
We assume that the input of a system consists of $n$ signals, one
for each participant. An input signal is denoted
$u=(u_1,\ldots,u_n)$, with $u_{i} \in \ell_{r_i,e}^{m_i}$ for some $m_i \in \mathbb N$ and
$r_i \in [1,\infty]$. 
A simple example is that of a dynamic system releasing at each period
the average over the past $l$ periods of the sum of the input values 
of the participants, i.e., with output
$\frac{1}{l} \sum_{k=t-l+1}^t \sum_{i=1}^n u_{i,k}$ at time $t$. 
For $r=(r_1,\ldots,r_n)$ and $m=(m_1,\ldots,m_n)$, an adjacency relation can be defined
on $l_{r,e}^m = \ell_{r_1,e}^{m_1} \times \ldots \times \ell_{r_n,e}^{m_n}$ by
$\Adj(u,u')$ if and only if $u$ and $u'$ differ by exactly one component signal, and moreover
this deviation is bounded. That is, let us fix a set of nonnegative numbers 
$b=(b_1, \ldots, b_n)$, $b_i \geq 0$, and define 
\begin{align}	\label{def: adjacency for vector signals}
\Adj^b(u,u') \text{ iff for some } i, \|u_i - u_i'\|_{r_i} \leq b_i, \\
\text{and } u_j = u_j' \text{ for all } j \neq i. \nonumber
\end{align}
Note that in (\ref{def: adjacency for vector signals}) two signals $u_i, u_i'$ are considered different 
if there exists some time $t$ at which $u_{i,t} \neq u'_{i,t}$.

\subsection{The Dynamic Gaussian Mechanism}

Recall (see, e.g., \cite{VanderSchaft00_passivity}) that for a system $F$ with inputs in $\ell_{r,e}^{m}$ 
and output in $\ell_{s,e}^{m'}$, its $\ell_{r}$-to-$\ell_{s}$ incremental gain $\gamma^{inc}_{r,s}(F)$ is defined 
as the smallest number $\gamma$ such that
\[
\| P_T F u- P_T F u' \|_s \leq \gamma \| P_T u - P_T u' \|_r, \;\; \forall u, u' \in \ell_{r,e}^m,
\]
for all $T$.
Now consider, for $r=(r_1,\ldots,r_n)$ and $m=(m_1,\ldots,m_n)$, a system $\mathcal G$ defined by
\begin{align}
&\mathcal G: l^m_{r,e} \to \ell_{s,e}^{m'} \nonumber \\
&\mathcal G (u_1,\ldots,u_n) = \sum_{i=1}^n \mathcal G_i u_i,	\label{eq: general system considered}
\end{align}
where $\mathcal G_i: \ell_{r_i,e}^{m_i} \to \ell_{s,e}^{m'}$, for all $1 \leq i \leq n$.
%

\begin{thm}	\label{thm: differentially private large scale system}
Let $\mathcal G$ be defined as in (\ref{eq: general system considered}) and
consider the adjacency relation (\ref{def: adjacency for vector signals}).
Then the mechanism $Mu = \mathcal Gu+w$, 
where $w$ is a white noise with $w_t \sim \mathcal N(0,\sigma^2 I_{m'})$
and $\sigma \geq \kappa(\delta,\epsilon) \max_{1 \leq i \leq n} \{ \gamma^{inc}_{r_i,2}(\mathcal G_i) \, b_i \}$,
is $(\epsilon,\delta)$-differentially private.
\end{thm}


 
 \begin{cor}	\label{cor: differentially private large scale system - linear}
 Let $\mathcal G$ be defined as in (\ref{eq: general system considered}) with each
 system $\mathcal G_i$ linear, and $r_i = 2$ for all $1 \leq i \leq n$.
Then the mechanism $Mu = \mathcal Gu+w$, where $w$ is a white Gaussian noise 
with $w_t \sim \mathcal N(0,\sigma^2 I_{m'})$ and
$\sigma \geq \kappa(\delta,\epsilon) \max_{1 \leq i \leq n} \{ \|\mathcal G_i\|_\infty \, b_i \}$,
is $(\epsilon,\delta)$-differentially private for  (\ref{def: adjacency for vector signals}).
\end{cor}

\subsection{Filter Approximation Set-ups for Differential Privacy}		\label{section: approximation set-ups}

\begin{figure}
\centering
\includegraphics[width=0.8\linewidth,height=2.5cm]{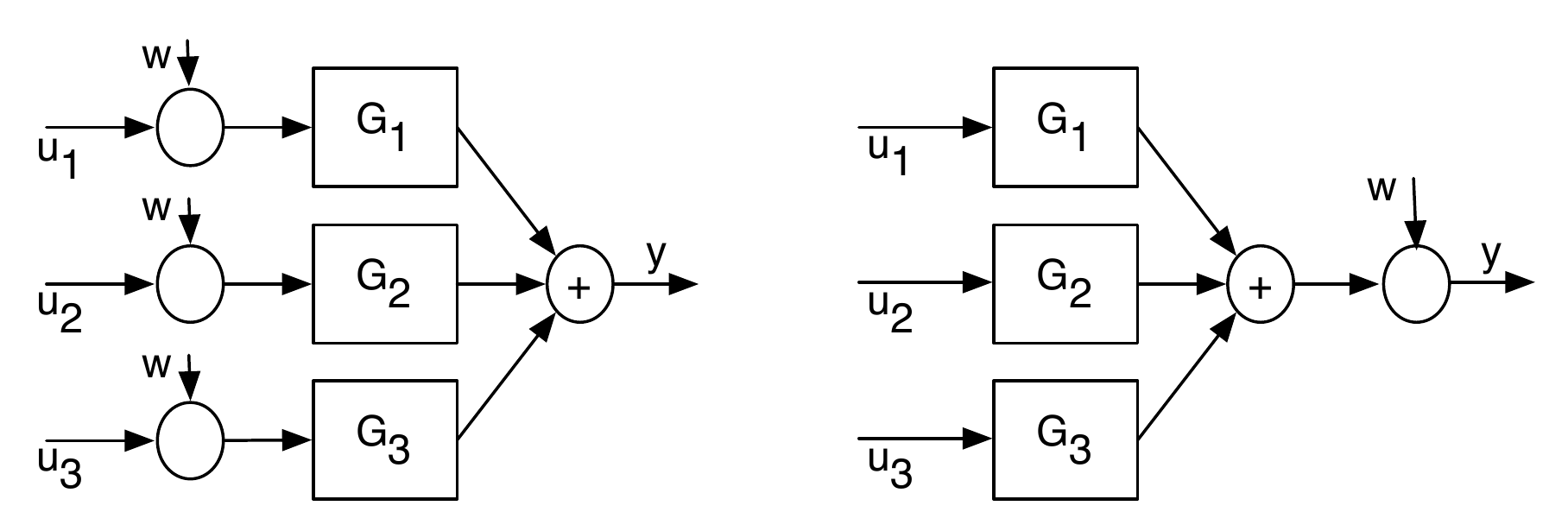}
\caption{Two architectures for differential privacy. (a) Input perturbation. (b) Output perturbation.}
\label{fig: DP basic architectures}
\end{figure}

Let $r_i = 2$ for all $i$ and $\mathcal G$ be linear as in the Corollary \ref{cor: differentially private large scale system - linear}, 
and assume for simplicity the same bound $b_1^2 = \ldots = b_n^2 = \mathcal E$
for the allowed variations in energy of each input signal. 
We have then two simple mechanisms producing a differentially private version of $\mathcal G$,
depicted on Fig.~\ref{fig: DP basic architectures}. 
The first one directly perturbs each input signal $u_i$ by adding to it a white Gaussian noise $w_i$
with $w_{i,t} \sim \mathcal N(0,\sigma^2 I_{m_i})$ and $\sigma^2 = \kappa(\delta,\epsilon)^2 \mathcal E$.
These perturbations on each input channel are then passed through $\mathcal G$,
leading to a mean squared error (MSE) for the output equal to 
$\kappa(\delta,\epsilon)^2 \mathcal E \|\mathcal G\|^2_2 
= \kappa(\delta,\epsilon)^2 \mathcal E \sum_{i=1}^n \|\mathcal G_i\|^2_2$. 
Alternatively, we can add a single source of noise at the output of $\mathcal G$ 
according to Corollary \ref{cor: differentially private large scale system - linear},
in which case the MSE is 
$\kappa(\delta,\epsilon)^2 \mathcal E \max_{1 \leq i \leq n} \{ \|\mathcal G_i\|_\infty^2 \}$.
Both of these schemes should be evaluated depending on the system $\mathcal G$
and the number $n$ of participants, as none of the error bound is better 
than the other in all circumstances. For example, if $n$ is small or if
the bandwidths of the individual transfer functions $\mathcal G_i$ do not overlap, 
the error bound for the input perturbation scheme can be smaller. 
Another advantage of this scheme is that the users can release 
differentially private signals themselves without relying on a trusted server.
However, there are cryptographic means for achieving the output perturbation
scheme without centralized trusted server as well, see, e.g., \cite{Shi11_DPaggregation}.


\begin{exmp}
Consider again the problem of releasing the average over the past $l$ periods
of the sum of the input signals, i.e., $\mathcal G = \sum_{i=1}^n \mathcal G_i$ with
\begin{align*}
(\mathcal G_i u_i)_t &= \frac{1}{l} \sum_{k=t-l+1}^t u_{i,k},
\end{align*}
for all $i$. Then $\|\mathcal G_i\|_2^2 = 1/l$, whereas $\|\mathcal G_i\|_\infty = 1$, for all $i$.
The MSE for the scheme with the noise at the input is then 
$\kappa(\delta,\epsilon)^2 \mathcal E n/l$. 
With the noise at the output, the MSE is $\kappa(\delta,\epsilon)^2 \mathcal E$,
which is better exactly when $n> l$, i.e., the number of users is larger than the
averaging window.
\end{exmp}

\section{Kalman Filtering}		\label{section: private KF}

We now discuss the Kalman filtering problem subject to a differential privacy constraint. 
With respect to the previous section, for Kalman filtering it is assumed that more is 
publicly known about the dynamics of the processes producing the individual signals. 
The goal here is to guarantee differential privacy for the individual state trajectories. 
Section \ref{section: traffic monitoring example} describes an application of the differentially private 
mechanisms presented here to a stylized traffic monitoring problem.

\subsection{A Differentially Private Kalman Filter}

Consider a set of $n$ linear systems, each with independent dynamics
\begin{align}	\label{eq: linear dynamics participant i}
x_{i,t+1} = A_i x_{i,t} + B_i w_{i,t}, \;\; t \geq 0, \;\; 1 \leq i \leq n,
\end{align}
where $w_i$ is a standard zero-mean Gaussian white noise process with covariance $\Exp[w_{i,t} w_{i,t'}] = \delta_{t-t'}$,
and the initial condition $x_{i,0}$ is a Gaussian random variable with mean $\bar x_{i,0}$, 
independent of the noise process $w_i$. 
System $i$, for $1 \leq i \leq n$, sends measurements 
\begin{align}	\label{eq: measurements participant i}
y_{i,t} = C_i x_{i,t} + D_i w_{i,t}
\end{align}
to a data aggregator. We assume for simplicity that the matrices $D_i$ are full row rank, 
and that $B_i D_i^T = 0$, i.e., the process and measurement noises are uncorrelated. 

\begin{figure}
\centering
\includegraphics[height=4cm]{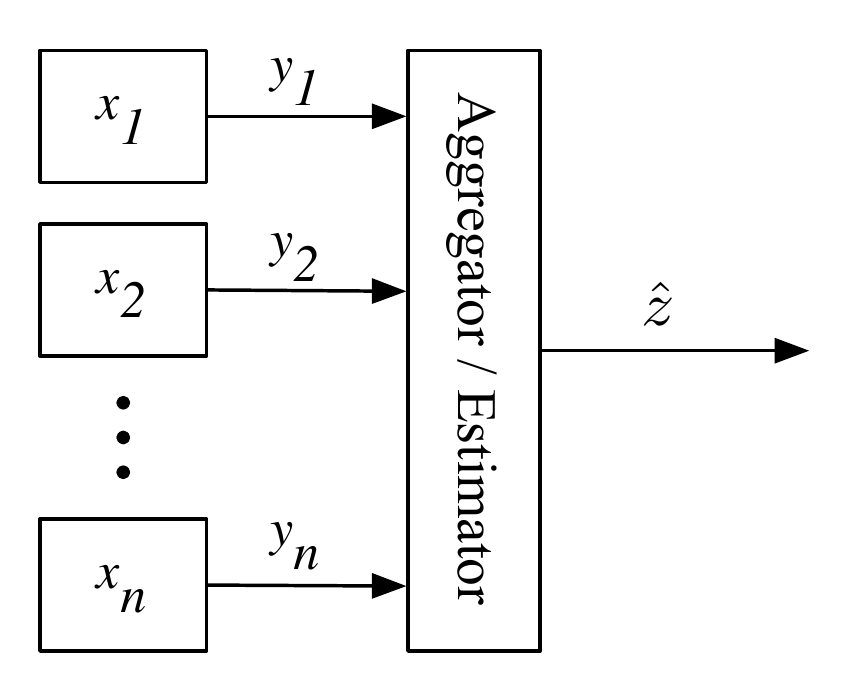}
\caption{Kalman filtering set-up.}
\label{fig: KF setup}
\end{figure}

%
The data aggregator aims at releasing a signal that asymptotically 
minimizes the minimum mean squared error with respect to a linear combination of the individual states. 
That is, the quantity of interest to be estimated at each period is $z_t = \sum_{i=1}^n L_i x_{i,t}$, 
where $L_i$ are given matrices, and we are looking for a causal estimator $\hat z$ constructed from the 
signals $y_i, 1 \leq i \leq n$, solution of
\[
\min_{\hat z} \lim_{T \to \infty} \frac{1}{T} \sum_{t=0}^{T-1} E \left[ \|z_t - \hat z_t\|_2^2 \right].
\]
The data $\bar x_{i,0}, A_i, B_i, C_i, D_i, L_i, 1 \leq i \leq n,$ are assumed to be public information.
For all $1 \leq i \leq n$, we assume that the pairs $(A_i, C_i)$ are detectable and the pairs
$(A_i, B_i)$ are stabilizable. 
%
%
In the absence of privacy constraint, the optimal estimator is $\hat z_t = \sum_{i=1}^n L_i \hat x_{i,t}$, 
with $\hat x_{i,t}$ provided by the steady-state Kalman filter \cite{Anderson05_filtering}.
Figure \ref{fig: KF setup} shows this initial set-up.

Suppose now that the publicly released estimate $\hat z$ should guarantee the differential privacy
of the participants. This requires that we first specify an adjacency relation on the appropriate space
of datasets. Let $x = [x_1^T, \ldots, x_n^T]^T$ and $y = [y_1^T, \ldots, y_n^T]^T$ denote
the global state and measurement signals. 
Assume that the mechanism is required to guarantee differential privacy with respect
to a subset $\mathcal S_i := \{i_1, \ldots, i_k\}$ of the coordinates of the state trajectory $x_i$.
Let the matrix $T_i$ be the diagonal matrix with $[T_i]_{jj} = 1$ if $j \in \mathcal S_i$, and $[T_i]_{jj} = 0$ otherwise. 
Hence $T_i v$ sets the coordinates of a vector $v$ which do not belong to the set $\{i_1,\ldots,i_k\}$ to zero.
Fix a vector $\rho \in \mathbb R_+^n$. The adjacency relation considered here is
\begin{align} \label{eq: adjacency for state trajectories - bis}
&\Adj_{\mathcal S}^\rho(x,x') \text{ iff } \text{for some } i, \; \|T_ix_i - T_ix_i' \|_2 \leq \rho_i, \\
& (I-T_i)x_i = (I-T_i)x_i', \text{and } x_j = x_j' \text{ for all } j \neq i. \nonumber
\end{align}
In words, two adjacent global output signals differ by the trajectory of a single participant, say $i$.
Moreover, for differential privacy guarantees we are constraining the range in energy variation
in the signal $T_i x_i$ of participant $i$ to be at most $\rho_i^2$. Hence, the distribution 
on the released results should be essentially the same if a participant's output signal value 
$T_i x_{i,t_0}$ at some single specific time $t_0$ were replaced by $T_i x'_{i,t_0}$ with 
$\|T_i (x_{i,t_0} - x'_{i,t_0})\| \leq \rho_i$, but the privacy guarantee should also hold for 
smaller instantaneous deviations on longer segments of trajectory. 

Depending on which signals on Fig. \ref{fig: KF setup} are actually published, and similarly to the discussion
of Section \ref{section: approximation set-ups}, there are different points at which a privacy inducing noise can be introduced.
First, for the input noise injection mechanism, the noise can be added by each participant directly 
to their transmitted measurement signal $y_i$.
Namely, since for two state trajectories $x_i, x_i'$ adjacent according to (\ref{eq: adjacency for state trajectories - bis})
we have for the corresponding measured signals
\[
\|y_i - y_i'\|_2 = \|C_i T_i (x_i - x_i') \|_2,
\]
differential privacy can be guaranteed if participant $i$ adds to $y_i$ a white Gaussian noise with covariance matrix 
$\kappa(\delta,\epsilon)^2 \rho_i^2 \sigma^2_{\max} (C_i T_i) I_{p_i}$, where $p_i$ is the dimension of $y_{i,t}$.
Note that in this sensitivity computation the measurement noise $D_i w_i$ has the same realization independently 
of the considered variation in $x_i$.
%
At the data aggregator, this additional noise can be taken into account in the design of the Kalman filter, 
since it can simply be viewed as an additional measurement noise.
%
Again, an important advantage of this mechanism is its simplicity of implementation when the participants 
do not trust the data aggregator, since the transmitted signals are already differentially private.

Next, consider the output noise injection mechanism. 
%
%
Since we assume that $\bar x^i_0$ is public information, the initial condition $\hat x_{i,0}$ of each state estimator is fixed. 
Consider now two state trajectories $x, x'$, adjacent according to (\ref{eq: adjacency for state trajectories - bis}), 
and let $\hat z, \hat z'$ be the corresponding estimates produced. We have
\[
\hat z - \hat z' = L_i \mathcal K_i (y_i - y_i') = L_i \mathcal K_i C_i T_i (x_i - x_i'),
\]
where $\mathcal K_i$ is the $i^{th}$ Kalman filter. Hence
$
\| \hat z - \hat z' \|_2 \leq \gamma_i \rho_i,
$
where $\gamma_i$ is the $\mathcal H_\infty$ norm of the transfer function $L_i \mathcal K_i C_i T_i$.
We thus have the following theorem.


\begin{thm}		\label{thm: general output perturbation mechanism result}
A mechanism releasing $\left( \sum_{i=1}^n L_i \mathcal K_i y_i \right) + \gamma \kappa(\delta,\epsilon) \; \nu$, 
where $\nu$ is a standard white Gaussian noise independent of $\{w_i\}_{1 \leq i \leq n}, \{x_{i,0}\}_{1 \leq i \leq n}$, 
and $\gamma = \max_{1 \leq i \leq n} \{\gamma_i \rho_i\}$, with $\gamma_i$ the $\mathcal H_\infty$ norm of 
$L_i \mathcal K_i C_i T_i$, is differentially private for the adjacency relation (\ref{eq: adjacency for state trajectories - bis}).
\end{thm}

\subsection{Filter Redesign for Stable Systems}	\label{section: filter redesign - stable}

In the case of the output perturbation mechanism, one can potentially improve the MSE
of the filter with respect to the Kalman filter considered in the previous subsection.
Namely, consider the design of $n$ filters of the form
\begin{align}
\hat x_{i,t+1} &= F_i \hat x_{i,t} + G_i y_{i,t}  	\label{eq: filter state} \\
\hat z_{i,t} &= H_i \hat x_{i,t} + K_i y_{i,t},	\label{eq: filter output}
\end{align}
for $1 \leq i \leq n$, where $F_i, G_i, H_i, K_i$ are matrices to determine. The estimator considered is
\[
\hat z_t = \sum_{i=1}^n \hat z_{i,t}, 
\]
so that each filter output $\hat z_{i}$ should minimize the steady-state error variance with $z_{i} = L_i x_{i}$,
and the released signal $\hat z$ should guarantee the differential privacy with respect to 
(\ref{eq: adjacency for state trajectories - bis}).
Assume first in this section that the system matrices $A_i$ are stable, in which case we also restrict the filter 
matrices $F_i$ to be stable.
Moreover, we only consider the design of full order filters, i.e., the dimensions of $F_i$ are greater or equal to those
of $A_i$, for all $1 \leq i \leq n$. Finally, we remove the simplifying assumption $B_i D_i^T = 0$.
%

Denote the overall state for each system and associated filter by $\tilde x_i = [x_i^T, \hat x_i^T]^T$.
The combined dynamics from $w_i$ to the estimation error $e_i := z_i - \hat z_i$ 
can then be written
\begin{align*}
\tilde x_{i,t+1} &=  \tilde A_i \tilde x_{i,t} + \tilde B_i w_{i,t} \\
e_{i,t} &= \tilde C_i \tilde x_{i,t} + \tilde D_i w_{i,t},
\end{align*}
where
\begin{align*}
\tilde A_i &= \begin{bmatrix}
A_i & 0 \\
G_iC_i & F_i
\end{bmatrix}, \;\;
\tilde B_i = \begin{bmatrix}
B_i \\ G_i D_i
\end{bmatrix}, \\
\tilde C_i &= \begin{bmatrix}
L_i - K_i C_i & -H_i
\end{bmatrix}, \;\;
\tilde D_i = -K_i D_i.
\end{align*}
The steady-state MSE for the $i^{th}$ estimator is then $\lim_{t \to \infty} \mathbb E[e_{i,t}^T e_{i,t}]$.

In addition, we are interested in designing filters with small $\mathcal H_\infty$ norm, in order to minimize the
amount of noise introduced by the privacy-preserving mechanism, which ultimately impacts the
overall MSE. Considering as in the previous subsection the sensitivity of filter $i$'s output to a change
from a state trajectory $x$ to an adjacent one $x'$ according to (\ref{eq: adjacency for state trajectories - bis}),
and letting $\delta x_i = x_i - x_i' = T_i (x_i - x_i') = T_i \delta x_i$, we see that the change in the output of filter $i$ follows the dynamics
\begin{align*}
\delta \hat x_{i,t+1} &= F_i \delta \hat x_{i,t} + G_i C_i T_i \delta x_i \\
\delta \hat z_i &= H_i \delta \hat x_{i,t} + K_i C_i T_i \delta x_i.
\end{align*}
Hence the $\ell_2$-sensitivity can be measured by the $\mathcal H_\infty$ norm of the transfer
function
\begin{align}	
\TF{F_i}{G_i C_i T_i}{H_i}{K_i C_i T_i}.
\end{align}

Simply replacing the Kalman filter in Theorem \ref{thm: general output perturbation mechanism result}, 
the MSE for the output perturbation mechanism guaranteeing $(\epsilon,\delta)$-privacy is then
\begin{align*}
\sum_{i=1}^n \|\tilde C_i (zI - \tilde A_i)^{-1} \tilde B_i + \tilde D_i \|_2^2 + \kappa(\delta,\epsilon)^2 
\max_{1\leq i \leq n} \{\gamma_i^2 \rho_i^2\}, \\
\text{with } \gamma_i := \| H_i (sI - F_i)^{-1} G_i C_i T_i + K_i C_i T_i \|_\infty.
\end{align*}
Hence minimizing this MSE leads us to the following optimization problem
\begin{align}
&\min_{\mu_i, \lambda, F_i, G_i, H_i, K_i} \quad 
\sum_{i=1}^n \mu_i + \kappa(\delta,\epsilon)^2 \lambda \label{eq: optimization filter} \\
& \text{s.t. } \forall \; 1 \leq i \leq n, \|\tilde C_i (zI - \tilde A_i)^{-1} \tilde B_i + \tilde D_i \|_2^2 \leq \mu_i, 
 \label{eq: optimization filter - constraint H2} \\
& \rho_i^2 \| H_i (zI - F_i)^{-1} G_i C_i T_i + K_i C_i T_i \|^2_\infty \leq \lambda.
\label{eq: optimization filter - constraint Hinf}
\end{align}
Assume without loss of generality that $\rho_i > 0$ for all $i$, since the privacy constraint 
for the signal $x_i$ vanishes if $\rho_i = 0$. 
The following theorem gives a convex sufficient condition in the form of Linear Matrix Inequalities (LMIs) 
guaranteeing that a choice of filter matrices $F_i, G_i, H_i, K_i$ satisfies the constraints 
(\ref{eq: optimization filter - constraint H2})-(\ref{eq: optimization filter - constraint Hinf}). These LMIs
can be obtained using the change of variable technique described in \cite{Scherer97_multiObj}.

\begin{thm}	\label{thm: LMI constraints filter design}
The constraints (\ref{eq: optimization filter - constraint H2})-(\ref{eq: optimization filter - constraint Hinf}), for some $1 \leq i \leq n$,
are satisfied if there exists matrices $W_i, Y_i, Z_i, \hat F_i, \hat G_i, \hat H_i, \hat K_i$ such that $\Tr(W_i) < \mu_i$, 
and the LMIs (\ref{eq: LMI1}), (\ref{eq: LMI2}) shown next page are satisfied.


\newcounter{MYtempeqncnt}
\begin{figure*}[!t]
\normalsize
\setcounter{MYtempeqncnt}{\value{equation}}
\setcounter{equation}{13}
\begin{align}	\label{eq: LMI1}
\begin{bmatrix}
W_i & (L_i - \hat K_i C_i - \hat H_i) & (L_i - \hat K_i C_i) & -\hat K_i D_i \\
* & Z_i & Z_i & 0 \\
* & * & Y_i & 0 \\
* & * & * & I
\end{bmatrix} \succ 0, \;\;
\begin{bmatrix}
Z_i & Z_i & 0 & 0 & 0 & 0 \\
* & Y_i & 0 & \hat F_i & 0 & \hat G_i C_i T_i \\
* & * & \frac{\lambda}{\rho^2_i} I & \hat H_i & 0 & \hat K_i C_i T_i \\
* & * & * & Z_i & Z_i & 0 \\
* & * & * & * & Y_i & 0 \\
* & * & * & * & * & I
\end{bmatrix} \succ 0,
\end{align}
\begin{align}	\label{eq: LMI2}
\begin{bmatrix}
Z_i & Z_i & Z_i A_i & Z_i A_i & Z_i B_i \\
* & Y_i & (Y_i A_i + \hat G_i C_i + \hat F_i) & (Y_i A_i + \hat G_i C_i) & (Y_i B_i + \hat G_i D_i) \\
* & * & Z_i & Z_i & 0 \\
* & * & * & Y_i & 0 \\
* & * & * & * & I
\end{bmatrix} \succ 0. 
\end{align}
%
\hrulefill
\vspace*{4pt}
\end{figure*}


If these conditions are satisfied, one can recover admissible filter matrices $F_i, G_i, H_i, K_i$ by setting
\begin{align}		\label{eq: recovering the filter}
F_i &= V_i^{-1} \hat F_i \hat Z_i^{-1} U_i^{-T}, \;\;
G_i = V_i^{-1} \hat G_i, \nonumber \\
H_i &= \hat H_i Z_i^{-1} U_i^{-T},
\;\; K_i = \hat K_i,
\end{align}
where $U_i, V_i$ are any two nonsingular matrices such that $V_i U_i^T = I - Y_i Z_i^{-1}$.
\end{thm}

Note that the problem (\ref{eq: optimization filter}) is also linear in $\mu_i, \lambda$. These variables can then
be minimized subject to the LMI constraints of Theorem \ref{thm: LMI constraints filter design} in order to design
a good filter trading off estimation error and $\ell^2$-sensitivity to minimize the overall MSE. 

\subsection{Unstable Systems}		\label{section: filter redesign - unstable}

If the dynamics (\ref{eq: linear dynamics participant i}) are not stable, the linear filter design approach
presented in the previous paragraph is not valid. To handle this case, we can further restrict the class of filters.
As before we minimize the estimation error variance together with the sensitivity measured 
by the $\mathcal H_\infty$ norm of the filter.
Starting from the general linear filter dynamics (\ref{eq: filter state}), (\ref{eq: filter output}), we can consider 
designs where $\hat x_i$ is an estimate of $x_i$, and set $H_i = L_i, K_i = 0,$ so that $\hat z_i = L_i \hat x_i$ is an estimate
of $z_i = L_i x_i$. The error dynamics $e_i := x_i - \hat x_i$ then satisfies
\[
e_{i,t+1} = (A_i-G_iC_i) x_{i,t} - F_i \hat x_{i,t} + (B_i - G_i D_i) w_{i,t}.
\]
Setting $F_i = (A_i - G_i C_i)$ gives an error dynamics independent of $x_i$
\begin{equation}	\label{eq: error dynamics - unstable case} 
e_{i,t+1} = (A_i-G_iC_i) e_{i,t} + (B_i - G_i D_i) w_{i,t},
\end{equation}
and leaves the matrix $G_i$ as the only remaining design variable. Note however that the resulting class of filters
contains the (one-step delayed) Kalman filter. To obtain a bounded error, there is an implicit constraint on $G_i$ that $A_i - G_i C_i$ 
should be stable. 

Now, following the discussion in the previous subsection, minimizing the MSE while enforcing differential privacy 
leads to the following optimization problem
\begin{align}
&\min_{\mu_i, \lambda, G_i} \quad \sum_{i=1}^n \mu_i + \kappa(\delta,\epsilon)^2 \lambda \label{eq: optimization filter - unstable} \\
&\text{s.t. } \; \forall \; 1 \leq i \leq n, \nonumber \\
& \| L_i (zI - (A_i - G_i C_i))^{-1} (B_i - G_i D_i) \|_2^2 \leq \mu_i, \label{eq: optimization filter - constraint H2 - unstable} \\
& \rho_i^2 \| L_i (zI - (A_i - G_i C_i))^{-1} G_i C_i T_i \|^2_\infty \leq \lambda. \label{eq: optimization filter - constraint Hinf - unstable}
\end{align}
Again, one can efficiently check a sufficient condition, in the form of the LMIs of the following theorem, 
guaranteeing that the constraints (\ref{eq: optimization filter - constraint H2 - unstable}), (\ref{eq: optimization filter - constraint Hinf - unstable}) 
are satisfied. Optimizing over the variables $\lambda_i, \mu_i, G_i$ can then be done using semidefinite programming. 

\begin{thm}	\label{thm: LMI constraints filter design - unstable}
The constraints (\ref{eq: optimization filter - constraint H2 - unstable})-(\ref{eq: optimization filter - constraint Hinf - unstable}), for some $1 \leq i \leq n$,
are satisfied if there exists matrices $Y_i, X_i, \hat G_i$ such that
\begin{align}	\label{H2 unstable LMI}
\Tr(Y_i L_i^T L_i) < \mu_i, \;\; 
\begin{bmatrix}
Y_i & I \\ I & X_i
\end{bmatrix} \succ 0, \\
\begin{bmatrix}
X_i & X_i A_i - \hat G_i C_i & X_i B_i - \hat G_i D_i \\
* & X_i & 0 \\
* & * & I
\end{bmatrix} \succ 0, 
\end{align}
\begin{align}			\label{Hinf unstable LMI}
\text{and }
\begin{bmatrix}
X_i & 0 & X_i A_i - \hat G_i C_i & \hat G_i C_i T_i  \\ 
* & \frac{\lambda}{\rho_i^2} I & L_i & 0 \\
* & * & X_i & 0 \\
* & * & * & I
\end{bmatrix} \succ 0.
\end{align}
If these conditions are satisfied, one can recover an admissible filter matrice $G_i$ by setting
\begin{align*} %
G_i = X_i^{-1} \hat G_i.
\end{align*}
\end{thm}


\section{A Traffic Monitoring Example}	\label{section: traffic monitoring example}

Consider a simplified description of a traffic monitoring system, inspired by real-world implementations
and associated privacy concerns as discussed in \cite{Sun04_trafficEstimation, Hoh11_VTL_trafficMonitoring} 
for example. There are $n$ participating vehicles traveling on a straight road segment. Vehicle $i$, for $1 \leq i \leq n$, is represented
by its state $x_{i,t} = [\xi_{i,t}, \dot \xi_{i,t}]^T$, with $\xi_i$ and $\dot \xi_i$ its position and velocity respectively.
This state evolves as a second-order system with unknown random acceleration inputs
\[
x_{i,t+1} = 
\begin{bmatrix}
1 & T_s \\ 
0 & 1
\end{bmatrix}
x_{i,t}
+ \sigma_{i1}
\begin{bmatrix}
T_s^2 / 2 & 0 \\ T_s & 0
\end{bmatrix}
w_{i,t}, 
\]
where $T_s$ is the sampling period, $w_{i,t}$ is a standard white Gaussian noise, and $\sigma_{i1} > 0$. 
Assume for simplicity that the noise signals $w_j$ for different vehicles are independent.
The traffic monitoring service collects GPS measurements from the 
vehicles \cite{Hoh11_VTL_trafficMonitoring}, thus getting noisy readings 
of the positions at the sampling times
\[
y_{i,t} = 
\begin{bmatrix}
1 & 0 
\end{bmatrix}
x_{i,t} 
+ 
\sigma_{i2}
\begin{bmatrix}
0 & 1 
\end{bmatrix}
w_{i,t},
\]
with $\sigma_{i2} > 0$.

The purpose of the traffic monitoring service is to continuously provide an estimate of the 
traffic flow velocity on the road segment, which is approximated by releasing 
at each sampling period an estimate of the average velocity of the participating vehicles, i.e., 
of the quantity
\begin{align}	\label{eq: sample average velocity}
z_t = \frac{1}{n} \sum_{i=1}^n \dot \xi_{i,t}.
\end{align}
With a larger number of participating vehicles, the sample average (\ref{eq: sample average velocity})
represents the traffic flow velocity more accurately. However, while individuals are generally interested 
in the aggregate information provided by such a system, e.g., to estimate their commute time, 
they do not wish their individual trajectories to be publicly revealed, since these might contain 
sensitive information about their driving behavior, frequently visited
locations, etc. 
The privacy mechanism proposed in \cite{Hoh11_VTL_trafficMonitoring} perturbs the GPS traces by
dropping $1$ out of $k$ measurements at each given location (the sampling is event based
rather than periodic as here). 
This makes individual trajectory tracking potentially harder, but no formal definition 
of privacy is introduced, and hence no quantitative privacy guarantee can be provided.

\subsection{Numerical Example} 

We now discuss some differentially private estimators introduced above, in the context of this example.
All individual systems are identical, hence we drop the subscript $i$ in the notation. Assume that the
selection matrix is  $T = \begin{bmatrix} 1 & 0 \\ 0 & 0 \end{bmatrix}$, that $\rho = 100$ m, $T_s=1 s$,
$\sigma_{i1} = \sigma_{i2} =  1$, and $\epsilon = \ln(3)$, $\delta = 0.05$. A single Kalman filter 
denoted $\mathcal K$ is designed to provide an estimate $\hat x_i$ of each state vector $x_i$, so that 
in absence of privacy constraint the final estimate would be 
\[
\hat z = \begin{bmatrix} 0 & \frac{1}{n} \end{bmatrix} \sum_{i=1}^n \mathcal K y_i = 
\begin{bmatrix} 0 & 1 \end{bmatrix} \mathcal K \left( \frac{1}{n} \sum_{i=1}^n y_i \right).
\]
Finally, assume that we have $n = 200$ participants, and that their mean initial velocity is $45$ km/h.

\begin{figure}
\centering
\includegraphics[height=6cm]{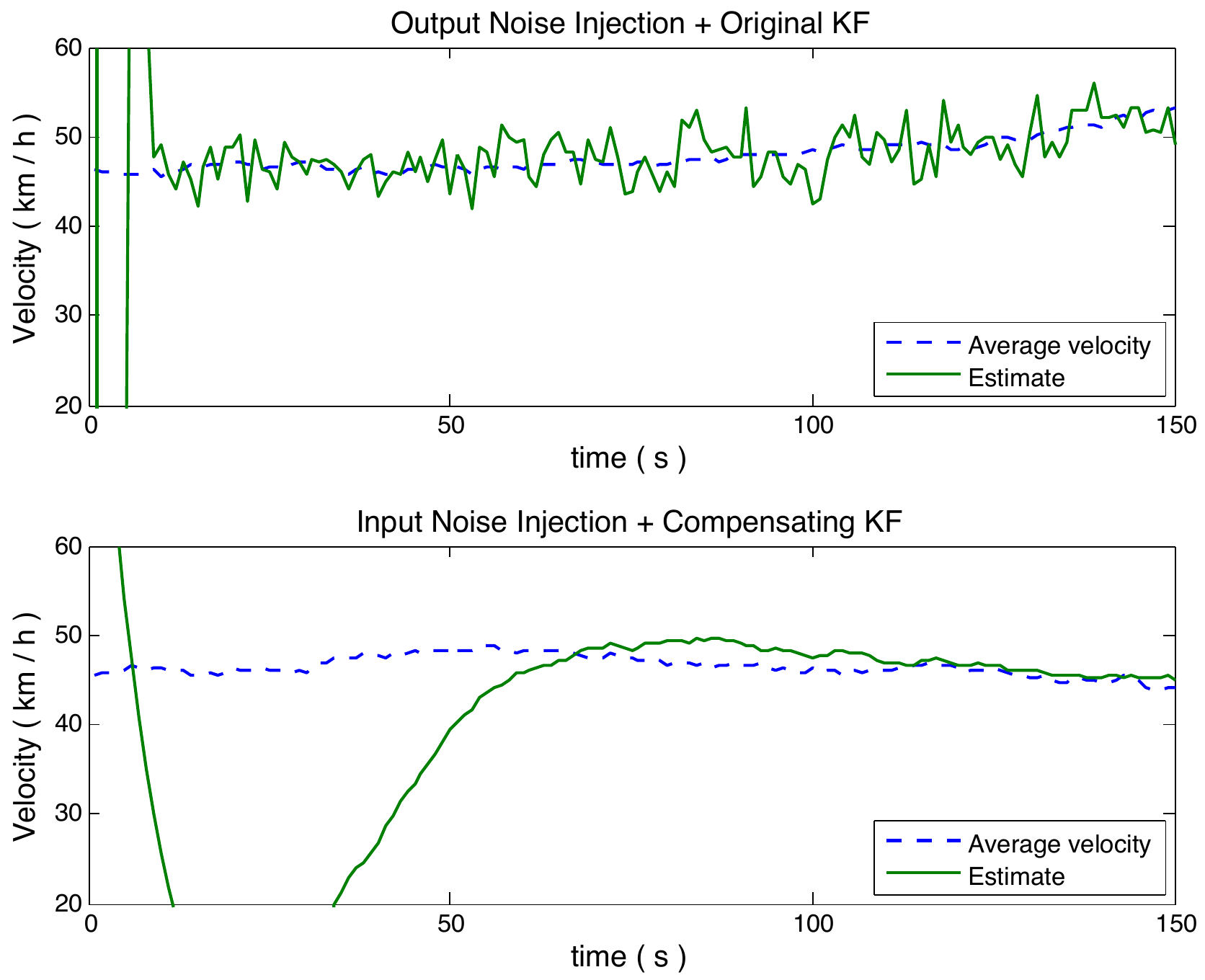}
\caption{Two differentially private average velocity estimates, with $n = 200$ users. 
The Kalman filters are initialized with the same incorrect initial mean velocity, in order to
evaluate their convergence time.}
\label{fig: velocity estimation}
\end{figure}

In this case, the input noise injection scheme without modification of the Kalman filter is essentially unusable
since its steady-state Root-Mean-Square-Error (RMSE) is almost $26$ km/h. However, modifying the Kalman filter to
take the privacy inducing noise into account as additional measurement noise leads to the best RMSE
of all the schemes discussed here, of about $0.31$ km/h.
Using the Kalman filter $\mathcal K$ with the output noise injection scheme leads to an RMSE of $2.41$ km/h. Moreover
in this case $\| \mathcal K \|_\infty = 0.57$ is quite small, and trying to balance estimation with sensitivity using the LMI
of Theorem \ref{thm: LMI constraints filter design - unstable} (by minimizing the MSE while constraining 
the $\mathcal H_\infty$ norm rather than using the objective function (\ref{eq: optimization filter - unstable}))
only allowed us to reduce this RMSE to $2.31$ km/h. 
However, an issue that is not captured in these steady-state estimation error measures is that of convergence time
of the filters. This is illustrated on Fig. \ref{fig: velocity estimation}, which shows a trajectory of the average 
velocity of the participants, together with the estimates produced by the input noise injection scheme with
compensating Kalman filter and the output noise injection scheme following $\mathcal K$. Although the RMSE
of the first scheme is much better, its convergence time of more than $1$ min, due to the large measurement noise assumed, 
is much larger. This can make this scheme impractical, e.g., if the system is supposed to respond quickly to an abrupt
change in average velocity.

\section{Conclusion}

We have discussed mechanisms for preserving the differential privacy of individual users
transmitting measurements of their state trajectories to a trusted central server releasing 
sanitized filtered outputs based on these measurements. Decentralized versions of these mechanisms 
can in fact be implemented in the absence of trusted server by means of cryptographic techniques \cite{Rastogi10_DPtimeSeries}.
Further research on privacy issues associated with emerging large-scale information processing and 
control systems is critical to encourage their development. 
Moreover, obtaining a better understanding of the design trade-offs between privacy or security 
and performance in these systems raises interesting system theoretic questions.






\bibliographystyle{IEEEtran}

\bibliography{IEEEabrv,/Users/jleny/Dropbox/Research/bibtex/energy,/Users/jleny/Dropbox/Research/bibtex/securityPrivacy,/Users/jleny/Dropbox/Research/bibtex/signalProcessing,/Users/jleny/Dropbox/Research/bibtex/controlSystems,/Users/jleny/Dropbox/Research/bibtex/communications}

\end{document}